\documentclass[10pt]{amsart}

\usepackage{amsthm,amsopn,verbatim}

\theoremstyle{plain}
\newtheorem*{mainthm}{Theorem \ref{MainResult}}

\newtheorem{thm}{Theorem}[section]
\newtheorem{lem}[thm]{Lemma}
\newtheorem{prop}[thm]{Proposition}
\newtheorem{cor}[thm]{Corollary}

\theoremstyle{definition}
\newtheorem{defn}{Definition}[section]

\theoremstyle{remark}
\newtheorem*{rem}{Remark}

\DeclareMathOperator{\Aut}{Aut}

\begin{document}

\title{Genericity of Filling Elements}
\author{Brent B. Solie}
\address{Department of Mathematics, University of Illinois at Urbana-Champaign, 1409 W. Green St., Urbana, IL 61801}
\email{solie@math.uiuc.edu}
\date{\today}
\subjclass[2010]{Primary 20P05, Secondary 03D15, 20E08, 20F65}
\keywords{geometric group theory, statistical group theory, free group, filling element, genericity}


\begin{abstract}

An element of a finitely generated non-Abelian free group $F(X)$ is said to be \emph{filling} if that  element has positive translation length in every very small minimal isometric action of $F(X)$ on an $\mathbb{R}$-tree.
We give a proof that the set of filling elements of $F(X)$ is exponentially $F(X)$-generic in the sense of Arzhantseva and Ol'shanski\u{\i}.
We also provide an algebraic sufficient condition for an element to be filling and show that there exists an exponentially $F(X)$-generic subset of filling elements whose membership problem is solvable in linear time.

\end{abstract} 

\maketitle

\section{Introduction}

Statistical group theory is a fairly new area of mathematics that has generated a great deal of active research in recent years.
A typical result in this area usually concerns the probability with which a random group element satisfies a given property.
Properties which are held by almost all sufficiently long elements of a given finitely generated group are especially interesting, and we informally refer to such properties as \emph{generic}.

The earliest occurence of the notion of genericity appears to be due to Guba \cite{Guba1986}.
Shortly afterwards, Gromov gave a formal definition in the context of finitely presented groups \cite{Gromov1987}.
In the same paper, Gromov asserts that almost every finitely presented group is hyperbolic, a fact first proved by Ol'shanski\u{\i} \cite{Olcprimeshanskiui1992} and later by Champetier \cite{Champetier1994,Champetier1995}.
Subsequent results in statistical group theory include the work of Arzhantseva \cite{Arzhantseva1998,Arzhantseva2000}, Arzhantseva and Ol'shanski\u{\i} \cite{Arzhantseva1996}, and Ollivier \cite{Ollivier2003,Ollivier2004,Ollivier2005a,Ollivier2005b}.
The surveys by Ghys \cite{Ghys2004} and Ollivier \cite{Ollivier2005} provide an excellent overview of genericity with a focus on random groups.
More recent results in statistical group theory apply the notion of genericity to computational group theory.
Some group-theoretic decision problems with high worst-case complexity have been shown to have low complexity on a generic set of inputs \cite{Kapovich2005a,Kapovich2006}.
These results have furthered the understanding of the average-case complexity of these problems \cite{Kapovich2003,Kapovich2005}.

The focus of the present paper is to investigate the generic geometric behavior of elements of a finitely generated non-Abelian free group.
Our investigation is motivated by the notion of a \emph{filling element}, first introduced by Kapovich and Lustig as a free group analogue of a filling curve on a surface \cite{Kapovich2010}.
Filling curves have been important in the theory of surface groups and may be characterized by a number of equivalent geometric and algebraic properties.

Let $\Sigma$ be a closed, orientable surface of genus at least two.
By a \emph{surface group} we mean the fundamental group of such a surface $\Sigma$.
Let $\alpha$ and $\beta$ be closed curves on $\Sigma$.
The geometric intersection number, denoted $i(\alpha,\beta)$, is the least number of intersections between members of the free homotopy classes of $\alpha$ and $\beta$.
If $\beta$ is such that $i(\alpha,\beta)>0$ for every essential simple closed curve $\alpha$, then we say that $\beta$ is a \emph{filling curve}.

Recall that the \emph{dual tree} associated to an essential simple closed curve $\alpha$ on $\Sigma$ is a simplicial tree equipped with a small minimal isometric action by $\pi_1(\Sigma)$.
(From now on, we will assume all of our group actions on trees are minimal and isometric.)
It is well-known that if $\beta$ is a (not necessarily simple) closed curve on $\Sigma$, then the \emph{translation length} of $\beta$ on $T_\alpha$, denoted $||\beta||_{T_\alpha}$, is equal to $i(\alpha,\beta)$.
Therefore, a closed curve $\beta$ is filling if and only if it has positive translation length on $T_\alpha$ for every essential simple closed curve $\alpha$.

As a consequence of Skora's duality theorem, any simplicial tree equipped with a small action by $\pi_1(\Sigma)$ can be collapsed down to a tree $T_\alpha$ for some essential simple closed curve $\alpha$.
Therefore, a closed curve $\beta$ is filling if and only if it has positive translation length in every small action of $\pi_1(\Sigma)$ on a simplicial tree.
An application of Bass-Serre theory shows that a closed curve $\beta$ is filling if and only if it is not conjugate into a vertex group of any elementary cyclic splitting of $\pi_1(\Sigma)$.

We now move from surface groups to consider a finitely generated non-Abelian free group $F(X)$.
Recall that the Culler-Vogtmann outer space, denoted $CV(F(X))$, is the projectivized space of free discrete actions of $F(X)$ on simplicial trees.
Outer space is the free group counterpart to Teichm\"{u}ller space in the sense that it admits a properly discontinuous isometric action by the outer automorphism group of $F(X)$ \cite{Culler1986}.
More, outer space also admits a Thurston-type compactification $\overline{CV}(F(X))$, where $\overline{CV}(F(X))$ is the projectivized space of very small actions of $F(X)$ on $\mathbb{R}$-trees \cite{Bestvina1994}.

In \cite{Kapovich2010}, Kapovich and Lustig introduce the notion of a filling element as a free group analogue for a filling curve.
A \emph{filling element} is an element $w \in F(X)$ that has positive translation length in every very small action of $F(X)$ on an $\mathbb{R}$-tree.
In the same paper, Kapovich and Lustig prove that almost every prefix of almost every freely reduced right-infinite word on the letters $X^\pm:=X \cup X^{-1}$ is a filling element of $F(X)$ \cite[Theorem 13.6]{Kapovich2010}.
This result echoes the work of Bonahon showing that filling is the typical behavior of closed curves on a closed orientable hyperbolic surface \cite{Bonahon1986,Bonahon1988,Bonahon1991}.

However, the proof of Kapovich and Lustig's theorem is non-constructive; we would still like to find an algorithmically verifiable sufficient condition for a free group element to be filling.
Further, we would like to find such a condition which is also generic in the following formal sense of Arzhantseva and Ol'shanski\u{\i}.
Let $S$ be the set of elements of $F(X)$ satisfying a given condition.
We say that $S$ is \emph{$F(X)$-generic} if 
\begin{displaymath}
  \displaystyle{\lim_{n \rightarrow \infty} \frac{\#(S \cap B_n)}{\# B_n}=1},
\end{displaymath}
where $B_n := \{ w \in F(X) : |w|_X \leq n \}$ \cite{Kapovich2006}.
If the limit converges to 1 exponentially fast, we say that $S$ is \emph{exponentially $F(X)$-generic}.
We will give a slightly more general definition of genericity in Definition \ref{Genericity}.

The main result of this paper is:

\begin{mainthm}
  Let $F(X)$ be a finitely generated non-Abelian free group.
  \begin{enumerate}
    \item
      Let $w \in F(X)$. If the stabilizer of $w$ in $\Aut F(X)$ is infinite cyclic, then $w$ is filling.
    \item
      The set of filling elements of $F(X)$ is exponentially $F(X)$-generic.
    \item
      There exists an exponentially $F(X)$-generic subset $S$ of $F(X)$ such that every element of $S$ is filling and the membership problem for $S$ is solvable in linear time.
  \end{enumerate}
\end{mainthm}

The proof of Theorem \ref{MainResult} relies on the construction of an exponentially $F(X)$-generic set $TS'$ constructed to study the generic-case complexity of Whitehead's algorithm \cite{Kapovich2006}.
Let $TS$ be the set of all cyclically reduced elements of $F(X)$ which are not proper powers, whose conjugacy classes are fixed by no relabeling Whitehead automorphism, and whose cyclic length is strictly increased by every non-inner, non-relabeling Whitehead automorphism.
(Whitehead automorphisms will be defined formally in Definition \ref{WHA}.)
The set $TS'$ is defined as the set of all elements of $F(X)$ whose cyclic reductions are elements of $TS$.
The set $TS'$ satisfies a number of important properties, outlined in Proposition \ref{TSProp}.

The other main ingredient in the proof of Theorem \ref{MainResult} is an analysis of the structure of vertex groups in elementary cyclic splittings of $F(X)$.
The work of Guirardel \cite{Guirardel1998} implies that a non-filling element must be conjugate into such a vertex group.
The algebraic structure of these vertex groups, detailed in in Proposition \ref{vertex}, implies part (1) of Theorem \ref{MainResult}.
Parts (2) and (3) then follow from part (1) and the properties of the set $TS'$.

\section{Preliminaries}

By $X$ we will always denote a finite set with $\# X \geq 2$. Let $F(X)$ be the free group on the set of letters $X$.

An \emph{$\mathbb{R}$-tree} is a geodesic metric space in which any two points are connected by a unique simple path.
We continue to assume that every action of $F(X)$ on an $\mathbb{R}$-tree is isometric and minimal.
An action of $F(X)$ on an $\mathbb{R}$-tree is \emph{very small} if the stabilizer of any tripod is trivial and the stabilizer of any arc is either trivial or maximal cyclic in the stabilizers of the endpoints of the arc.
For $w \in F(X)$ and $T$ an $\mathbb{R}$-tree on which $F(X)$ acts, then the \emph{translation length} of $w$ is defined to be
\begin{displaymath}
  ||w||_T := \displaystyle{\inf_{p \in T} d_T(p,w(p))}.
\end{displaymath}
\begin{defn}
  Let $w \in F(X)$. We say that $w$ is a \emph{filling element} if $||w||_T > 0$ for every very small action of $F(X)$ on an $\mathbb{R}$-tree $T$. We say that $w$ is a \emph{non-filling element} if $||w||_T = 0$ for some very small action of $F(X)$ on an $\mathbb{R}$-tree $T$.
\end{defn}

The work of Guirardel allows us to approximate the very small action of $F(X)$ on a given $\mathbb{R}$-tree by very small actions on a simplicial trees \cite{Guirardel1998}.
In particular, if $w \in F(X)$ fixes a point in a very small action on an $\mathbb{R}$-tree, it must also fix a point in a very small action on a simplicial tree.
This implies:

\begin{prop}
  \label{Guir}
  An element $w \in F(X)$ is non-filling if and only if $||w||_T = 0$ for some very small action of $F(X)$ on a simplicial tree $T$.
\end{prop}

Recall that a very small action of $F(X)$ on a simplicial tree gives a decomposition of $F(X)$ as the fundamental group of a graph of groups with cyclic edge groups. We briefly review some of the associated terminology.

\begin{defn}
  A \emph{cyclic splitting} of $F(X)$ is the decomposition of $F(X)$ as the fundamental group of a graph of groups with cyclic edge groups.
  A \emph{boundary map} in a cyclic splitting is a map from an edge group to a vertex group.
  A cyclic splitting of $F(X)$ is \emph{elementary} if the corresponding graph of groups is connected and has exactly one edge.
  An elementary cyclic splitting of $F(X)$ is a \emph{segment of groups} if it has two distinct vertices and is a \emph{loop of groups} if it has a single vertex.
  An elementary cyclic splitting is \emph{nontrivial} if it is either a loop of groups, or it is a segment of groups in which neither boundary map is an isomorphism.
  Given any splitting of $F(X)$, we say that $w \in F(X)$ is \emph{elliptic} with respect to this splitting if $w$ is conjugate to an element of a vertex group.
  Elements of $F(X)$ which are not elliptic in a given splitting are said to be \emph{hyperbolic}.
\end{defn}

We will find it useful to slightly generalize the notion of $F(X)$-genericity presented in the introduction. The following notion of genericity referred to as the \emph{Arzhantseva-Ol'shanski\u{\i} model of genericity}, and the specific terminology is due to Kapovich, Schupp, and Shpilrain \cite{Kapovich2006}.

\begin{defn}
\label{Genericity}
Let $S \subseteq T \subseteq F(X)$.
We say that $S$ is \emph{$T$-generic} if
\begin{displaymath}
  \displaystyle{\lim_{n \rightarrow \infty} \dfrac{\#(S \cap B_n)}{\#(T \cap B_n)} = 1},
\end{displaymath}
where $B_n := \{w \in F(X) : |w|_X \leq n\}$ and $|w|_X$ denotes the word length of the element $w$ with respect to the basis $X$.
If the above limit converges exponentially fast, we say that $S$ is \emph{exponentially $T$-generic}.
If $T-S$ is (exponentially) $T$-generic, then we say that $S$ is \emph{(exponentially) $T$-negligable}.
\end{defn}

Recall that for the standard basis $X$ of $F(X)$, we have a corresponding set of Whitehead automorphisms which act as a finite generating set for $\Aut F(X)$.
Each of these Whitehead automorphism falls into one of the two following categories.

\begin{defn}
\label{WHA}
  An automorphism $\sigma: F(X) \rightarrow F(X)$ is called a \emph{type I Whitehead automorphism} (or a \emph{relabeling Whitehead automorphism}) if $\sigma$ is induced by a permutation on the set $X^\pm := X \cup X^{-1}$.

  An automorphism $\sigma: F(X) \rightarrow F(X)$ is called a \emph{type II Whitehead automorphism} if there is an element $a \in X^\pm$ such that for all $x \in X^\pm$, we have $\sigma(x) \in \{x, xa, a^{-1}x, x^a \}$, where $x^a := z^{-1}xa$.
\end{defn} 

\section{Main Results}

\subsection{Cyclic Splittings of $F(X)$}

The structure of elementary cyclic splittings of free groups has been well-studied.
A topological lemma due to Bestvina and Feighn \cite[Lemma 4.1]{Bestvina1994} can be used to characterize one-relator presentations of $F(X)$, resulting in the following characterization of vertex groups of cyclic splittings.
The following proposition may also be obtained from earlier results in \cite{Shenitzer1955}, \cite{Stallings1980}, or \cite{Swarup1986}.

\begin{prop}
  \label{vertex}
  If $F(X)$ has the structure of a nontrivial segment of groups with cyclic edge group, then its vertex groups must have the form $\langle A, b \rangle$ and $\langle B \rangle$, where $A \sqcup B$ is a basis for $F(X)$, $\#B \geq 2$, and $b \in \langle B \rangle$.

  If $F(X)$ has the structure of a nontrivial loop of groups with cyclic edge group, then its vertex group must have the form $\langle U, u^v \rangle$, where $U \sqcup \{v\}$ is a basis for $F(X)$ and $u \in \langle U \rangle$.

\end{prop}

\subsection{The Set $TS'$}

In \cite{Kapovich2006}, Kapovich, Schupp, and Shpilrain construct an exponentially $F(X)$-generic set with several important properties related to Whitehead's algorithm.

\begin{defn}
  Let $C \subseteq F(X)$ be the set of cyclically and freely reduced elements of $F(X)$.
  The set \emph{$TS$} is the set of $w \in C$ which are not proper powers, whose cyclic length is increased by every non-inner type II Whitehead automorphism, and whose conjugacy class is fixed by no type I Whitehead automorphism.
  The set $TS'$ is the set of elements $w \in F(X)$ whose cyclic reductions are in $TS$.
\end{defn}

\begin{prop}[{\cite[Theorem 8.5]{Kapovich2006}}]
  \label{TSProp}
  Let $\#X \geq 2$ and let $TS' \subseteq F(X)$ be as above.
  \begin{enumerate}
    \item
      The set $TS'$ is exponentially $F(X)$-generic.
    \item
      For any nontrivial $w \in TS'$, the stabilizer of $w$ in $\Aut F(X)$ is the infinite cyclic group generated by right-conjugation by $w$.
    \item
      The membership problem for $TS'$ is solvable in linear time.
  \end{enumerate}
\end{prop}

  We very briefly summarize the arguments presented in \cite{Kapovich2006}.
  Let $N:=\#X$.
  Let $w \in C$ have cyclic length $n$.
  For every $x \in X^\pm$, define $w_x$ to be the number of occurrences of $x$ in $w$.
  Likewise, for every $x, y \in X^\pm$ such that $x \neq y^{-1}$, define $w_{xy}$ to be the number of occurrences of $xy$ as a subword of $w$ written as a cyclic word.

  Let $\epsilon > 0$. Let $L(\epsilon)$ be the set of all $w \in C$ such that for all $x,y \in X^\pm$ where $x \neq y^{-1}$, we have
  \begin{align*}
    \dfrac{w_x}{n} \in \left(\dfrac{1}{2N} - \epsilon, \dfrac{1}{2N} + \epsilon \right)  \text{\: and \:}
    \dfrac{w_{xy}}{n} \in \left(\dfrac{1}{2N(2N-1)} - \epsilon, \dfrac{1}{2N(2N-1)} + \epsilon \right).
  \end{align*}

  If $\tau$ is a type II Whitehead automorphism, the difference in cyclic length between $w$ and $\tau(w)$ is easily calculated in terms of the quantities $w_x$ and $w_{xy}$ (see \cite{Lyndon2001}, \cite{Roig2007}).
  Let $0< \epsilon < (2N-3)/N(2N-1)(4N-3)$.
  If $w \in L(\epsilon)$ for such an $\epsilon$, the quantities $w_x$ and $w_{xy}$ nearly uniformly distributed.
  It is then fairly straightforward to see that any non-inner type II Whitehead automorphism will necessarily increase the cyclic length of $w \in L(\epsilon)$.

  Using some deep results from large deviation theory, one can further show that the set $L(\epsilon)$ itself is $C$-generic.
  Since $L(\epsilon)$ is $C$-generic, the set $L'(\epsilon)$ consisting of elements of $F(X)$ whose cyclic reductions are in $L(\epsilon)$ is $F(X)$-generic.
  Although $L'(\epsilon)$ may contain elements whose conjugacy class is fixed by a type I Whitehead automorphism, the set of such elements is $F(X)$-negligable and may be safely discarded to obtain a subset of $TS'$ which is $F(X)$-generic.
  We conclude that $TS'$ itself is also $F(X)$-generic.

  Suppose now that $w \in TS$.
  If $\alpha \in \Aut F(X)$ is such that $\alpha(w) = w$, recall that Whitehead's theorem states that $\alpha$ can be written as a sequence of Whitehead automorphisms $\alpha = \tau_m \circ \dots \circ \tau_1$ such that $\tau_{i+1}$ does not increase the cyclic length of $\tau_i \circ \dots \tau_1 (w)$.
  Since $w \in TS$ and $TS$ is closed under type I Whitehead automorphisms and cyclic permutations, each $\tau_i$ must be either a type I Whitehead automorphism or it must cyclically permute $w$.
  We may therefore write $\alpha = \gamma \circ \sigma$, where  $\sigma$ is a type I Whitehead automorphism and $\gamma$ is an inner automorphism.
  Since the conjugacy class of $w$ cannot be fixed by a type I Whitehead automorphism, $\sigma$ must be trivial, and so $\alpha$ itself must be inner.
  It follows directly that $\alpha$ can only be conjugation by some power of $w$, since $w$ is not a proper power itself.
  A slight modification to this argument shows that any $w'$ whose cyclic reduction is $w$ also has a cyclic stabilizer in $\Aut F(X)$.

  To algorithmically decide whether an element $w'$ belongs to $TS'$, we first pass to the cyclic reduction of $w'$, denoted by $w$, in time linear in the length of $w'$.
  We then decide whether $w \in TS$.
  It is well-known how to check whether $w$ is a proper power or not in time linear in the length of $w$.
  Checking whether the cyclic length of $w$ increases under any non-inner type II Whitehead automorphisms or whether any type I Whitehead automorphism fixes the conjugacy class of $w$ can each be done in linear time since the number of Whitehead automorphisms of $F(X)$ is constant.

\begin{rem}
  Although the membership problem for $TS'$ is linear in the length of the input, it has a factorial dependence on the rank of $F(X)$ due to the growth of the set of Whitehead automorphisms.
  It is worth noting that we can circumvent this dependence by passing to $L'(\epsilon)$, which is exponentially $TS'$-generic for the sufficiently small values of $\epsilon$ given above.
  The membership problem for $L'(\epsilon)$ requires only that we check the counts on a given word's one- and two-letter subwords, which can be done in linear time without any dependence on the rank of the ambient free group.
\end{rem}

\subsection{Genericity of Filling Elements}

We first consider the case where $w \in F(X)$ is elliptic in an elementary splitting of $F(X)$ over the trivial group.

\begin{lem}
  \label{FFStab}
  Let $w \in F(X)$ be elliptic in an elementary splitting of $F(X)$ over a trivial group.
  Then $w$ has non-cyclic stabilizer in $\Aut F(X)$.
\end{lem}

\begin{proof}
  To say that $w$ is elliptic in an elementary splitting of $F(X)$ over a trivial group is equivalent to saying that $w$ is contained in a proper free factor of $F(X)$.
  Suppose that $w$ is not a proper power.
  Let $A \sqcup B$ be a basis for $F(X)$ such that $\#A, \#B \geq 1$ and $w \in \langle A \rangle$.
  Let $\sigma : F(X) \rightarrow F(X)$ be right-conjugation by $w$.
  Define $\tau: F(X) \rightarrow F(X)$ via
  
  \begin{equation*}
    \tau(x) =
      \begin{cases}
        x^w & \text{if } x \in A, \\
        x & \text{if } x \in B,
      \end{cases}
  \end{equation*}

  where $x^w := w^{-1}xw$.
  Both $\sigma$ and $\tau$ fix $w$. However, $\sigma$ fixes exactly $\langle w \rangle$, while $\tau$ fixes $\langle w, B \rangle$.
  Thus $\sigma$ must be distinct from every power of $\tau$, so the $\Aut F(X)$ stabilizer of $w$ cannot be cyclic.

  If $w=z^r$ where $r > 1$ and $z$ is not a proper power, then $z$ is elliptic in an elementary cyclic splitting if and only if $w$ is elliptic in that same splitting.
  We may therefore pass from $w$ to its root $z$, which is also elliptic in the given splitting.
  The argument above shows that $z$ has a non-cyclic stabilizer, and since the stabilizer of $w$ contains that of $z$, $w$ must have non-cyclic stabilizer in $\Aut F(X)$ as well.
\end{proof}

Since the set $TS'$ is an exponentially $F(X)$-generic set whose elements all have cyclic stabilizers in $\Aut F(X)$, any set consisting of elements with non-cyclic stabilizers is exponentially $F(X)$-negligable.

\begin{cor}
  \label{FFNeg}
  The set of elements of $F(X)$ which lie in a proper free factor of $F(X)$ is exponentially $F(X)$-negligable.
\end{cor}

\begin{rem}
  This is a slight generalization of results appearing in \cite{Borovik2002} and \cite{Burillo2002}, which show that the set of primitive elements of $F(X)$ is $F(X)$-negligable.
\end{rem}

\begin{lem}
  \label{Stab}
  Let $w \in F(X)$ be elliptic in some elementary cyclic splitting of $F(X)$.
  Then $w$ has a non-cyclic stabilizer in $\Aut F(X)$.
\end{lem}

\begin{proof}
  Suppose that $w$ is not a proper power.

  Let $w \in F(X)$ be elliptic in a segment of groups.
  Then there must exist a basis $A \sqcup B$ of $F(X)$ such that $\#A \geq 1$, $\#B \geq 2$, $b \in \langle B \rangle$, and either $w \in \langle A, b \rangle$ or $w \in \langle B \rangle$.
  Note that if $b$ is a proper power of some $c \in F(X)$, then we would have $w \in \langle A, c \rangle$, so $w$ would remain elliptic in a splitting of the same type.
  Hence we may assume that $b$ is not a proper power.
  Define an automorphism $\sigma: F(X) \rightarrow F(X)$ by
  \begin{displaymath}
    \sigma(y) =
    \begin{cases}
      y, & \text{if } y \in A \\
      y^b, & \text{if } y \in B.
    \end{cases}
  \end{displaymath}

  Any power of $\sigma$ fixes the rank 2 subgroup $\langle A, b \rangle$ pointwise and so also fixes $w$, whereas right-conjugation by $w$ fixes exactly the cyclic subgroup $\langle w \rangle$.
  Right-conjugation by $w$ must therefore differ from every power of $\sigma$, so the stabilizer of $w$ in $\Aut F(X)$ cannot be cyclic.

  If $w \in \langle B \rangle$, since $\#A \geq 1$, $w$ lies in a proper free factor of $F(X)$.
  Lemma \ref{FFStab} states that such an element has a non-cyclic stabilizer in $\Aut F(X)$.

  Let $w \in F(X)$ be elliptic in a loop of groups.
  There then exists a basis $U \sqcup \{v\}$ of $F(X)$ such that $w \in \langle U, u^v \rangle$ for some $u \in \langle U \rangle$.
  We define the map $\tau: F(X) \rightarrow F(X)$ by
  \begin{align*}
    \tau(y) &= y, \text{  if } y \in U \\
    \tau(v) &= uv.
  \end{align*}
  
  Since $u \in \langle U \rangle$, $\tau$ is an automorphism.
  In particular, $\tau$ fixes the subgroup $\langle U, u^v \rangle$ pointwise, so no power of $\tau$ equals right-conjugation by $x$, which fixes only the cyclic subgroup $\langle w \rangle$.
  Again, the stabilizer of $w$ in $\Aut F(X)$ therefore cannot be cyclic.

  We handle the case where $w$ is a proper power in the same way it was handled in the proof of Lemma \ref{FFStab}.
\end{proof}

\begin{lem}
  If $w \in TS'$, then $w$ is filling.
\end{lem}

\begin{proof}
  We consider instead the set of non-filling elements in $F(X)$.
  By Proposition \ref{Guir}, a non-filling element $w$ fixes a point in some simplicial tree equipped with a very small action by $F(X)$.
  This in turn gives a splitting of $F(X)$ over infinite cyclic and trivial groups.
  Since $w$ fixes a point in the action, $w$ must be elliptic in this splitting.
  We may choose any edge in the splitting and collapse the components of its complement down to vertices, thereby obtaining an elementary splitting of $F(X)$ over an infinite cyclic or trivial group in which $w$ is elliptic.
  Since $w$ is elliptic in such a splitting, by Lemmas \ref{FFStab} and \ref{Stab}, the stabilizer of $w$ is non-cyclic in $\Aut F(X)$.
  By Proposition \ref{TSProp}, all elements of $TS'$ have cyclic stabilizers in $\Aut F(X)$, and so $w \notin TS'$.
\end{proof}

\begin{thm}
  \label{MainResult}
  Let $F(X)$ be a finitely generated non-Abelian free group.
  \begin{enumerate}
    \item
      Let $w \in F(X)$. If the stabilizer of $w$ in $\Aut F(X)$ is infinite cyclic, then $w$ is filling.
    \item
      The set of filling elements of $F(X)$ is exponentially $F(X)$-generic.
    \item
      There exists an exponentially $F(X)$-generic subset $S$ of $F(X)$ such that every element of $S$ is filling and the membership problem for $S$ is solvable in linear time.
  \end{enumerate}
\end{thm}

\begin{proof}
  Part (1) follows from Lemmas \ref{FFStab} and \ref{Stab}.
  Since every element of $TS'$ has a cyclic stabilizer in $\Aut F(X)$ (Proposition \ref{TSProp}, part (1)), every element of $TS'$ must be filling.
  Part (2) then follows from the fact that $TS'$ is exponentially $F(X)$-generic (Proposition \ref{TSProp}, part (2)).
  Finally, Part (3) follows from part (3) of Proposition \ref{TSProp}, taking $S$ to be exactly $TS'$.
\end{proof} 

\section{Acknowledgements}

The author would like to thank Ilya Kapovich for his guidance and many valuable comments and Enric Ventura for a helpful discussion via email. 

\bibliographystyle{plain}
\bibliography{genbib}

\end{document}